  \documentclass[12pt]{amsart}
  \usepackage{amssymb,amsthm}
\theoremstyle{plain}
  \newtheorem{theorem}{Theorem}[section]
  \newtheorem{lemma}[theorem]{Lemma}
  \newtheorem{corollary}[theorem]{Corollary}
  \newtheorem{proposition}[theorem]{Proposition}
  \newtheorem*{theoA}{Theorem A}
  \newtheorem*{theoB}{Theorem B}

\theoremstyle{definition}
  \newtheorem{remark}[theorem]{Remark}

     \renewcommand{\theenumi}{({\thetheorem}.\arabic{enumi})}
  \newenvironment{numpar}
    {\stepcounter{theorem}
     \smallskip
     \par {\bf ({\thetheorem})}}

 \newenvironment{remarks}
   {\stepcounter{theorem} \smallskip   \noindent  {\bf Remarks.}
      \begin{list}{\bf {\theenumi}.}
       {\setlength{\leftmargin}{0in}
      \setlength{\rightmargin}{0in} \setlength{\topsep}{0in}
      \setlength{\parsep}{0in} \setlength{\itemindent}{6em}}
      \usecounter{enumi} }
   {\end{list}          \medskip}

 \newenvironment{examples}
   {\stepcounter{theorem} \smallskip   \noindent  {\bf Examples.}
      \begin{list}{\bf {\theenumi}.}
       {\setlength{\leftmargin}{0in}
      \setlength{\rightmargin}{0in} \setlength{\topsep}{0in}
      \setlength{\parsep}{0in} \setlength{\itemindent}{6em}}
      \usecounter{enumi} }
   {\end{list}          \medskip}

  \newcommand{\Ical}{\mathcal{I}}
  \newcommand{\Mcal}{\mathcal{M}}
  \newcommand{\Ocal}{\mathcal{O}}
  \newcommand{\Cbb}{\mathbb{C}}
  \newcommand{\Pbb}{\mathbb{P}}
  \newcommand{\Fbb}{\mathbb{F}}

  \newcommand{\codim}{\operatorname{codim}}
  \newcommand{\corank}{\operatorname{corank}}
  \newcommand{\crk}{\operatorname{crk}}
  \newcommand{\Hom}{\operatorname{Hom}}
  \newcommand{\Ker}{\operatorname{Ker}}
  \newcommand{\Supp}{\operatorname{Supp}}

\begin{document}

\title{Hyperplane Sections of Calabi-Yau Varieties}
\author{Jonathan Wahl}
\address{Department of Mathematics \\
         University of North Carolina \\
         Chapel Hill, NC 27599--3250   }
\email{jw@math.unc.edu}
\subjclass{AMS Subject Classification Numbers
14J32, 14D15}
\dedicatory{In Memory of My Friend Michael Schneider}

  \begin{abstract}
  \smallskip

 {\bf Theorem}. \textit{If $W$ is a smooth complex projective
variety with $h^1(\Ocal_W)=0$, then a sufficiently ample smooth
divisor $X$ on $W$ cannot be a hyperplane section of a Calabi-Yau
variety, unless $W$ is itself a Calabi-Yau}.

 {\bf Corollary}. \textit{ A smooth hypersurface of degree $d$ in
$\mathbb{P}^n$ ($n\geq 2$) is a hyperplane section of a Calabi-Yau
variety iff $n+2\leq d\leq 2n+2$}.

The method is to construct out of the variety $W$ a universal family
of all varieties $Z$ for which $X$ is a hyperplane section with
normal bundle $K_X$, and examine the ``bad'' singularities of such
$Z$.

It was proved in [W1] that if a smooth curve lies on a $K$-3
surface, its Gaussian-Wahl map $\Phi_K$ is not surjective.

 {\bf Theorem}. \textit{The following smooth curves do not lie on a
$K$-3,
even though $\Phi_K$ is not surjective: plane curves of degree $\geq 7$;
bielliptic curves of genus $\geq 11$; curves on $\mathbb{F}_n$ of
degree $\geq 5$ over $\mathbb{P}^1$}.
  \end{abstract}

 \maketitle

  \newpage
  \section*{0. Introduction}

A {\em Calabi-Yau variety} shall mean a complex projective variety
$Y$ with trivial dualizing sheaf $K_Y \cong \Ocal_Y$ and
$h^1(\Ocal_Y)=0$, with only isolated canonical (i.e. rational)
singularities.  When $\dim Y=2$, $Y$ is a $K$-3 surface, possibly
with rational double points.  If $X$ is a smooth hyperplane section
of a Calabi-Yau, by the adjunction formula $K_X$ is very ample; if
$\dim Y\geq 3$ also $h^1(\Ocal_X)=0$.

     Conversely, given a canonically polarized $(X,K_X)$, one asks if
it can be a hyperplane section of a Calabi-Yau.  Expected finiteness
results for families of Calabi-Yaus suggest such $(X,K_X)$ should be
quite special.  For instance, there is a 19-dimensional family of
$K$-3 surfaces with a hyperplane section of genus $g$, hence only a
$(g+19)$-dimensional family of such curves [MoM].  Still, $K$-3
curves can be Brill-Noether-Petri general [L].

     By considering deformations of the affine cone $A$ over a
canonical curve, we showed in [W1] that a curve $X$ on a $K$-3 has
non-surjective Gaussian-Wahl map $\Phi_K: \Lambda^2(K)\rightarrow
\Gamma (K^{\otimes 3})$; the relevant part of the tangent space to
the deformations of $A$ must be non-$0$. Curves with surjective
$\Phi_K$ include most complete intersections and the generic curve
of genus 10 or $\geq 12$ [CHM]; such $X$ could not be on a $K$-3.
Conversely, it is conjectured in [W2] that a Brill-Noether-Petri
general curve with $\Phi_K$ non-surjective does indeed lie on a
$K$-3.

        This paper shows many $X$ are {\em not} hyperplane sections
of a Calabi-Yau, even when there are ``interesting'' deformations
of the cone. For instance, smooth plane curves of degree $\geq 7$
cannot lie on a $K$-3 [GL], although $\Phi_K$ has corank 10.  For
$\dim X\geq 2$, the relevant tangent space is
$H^1(X,\Theta_X\otimes K_X^{-1})$, which is always non-0 for
$X\subset \Pbb^n$ a smooth hypersurface; but we prove that for
degree $>2n+2$, such an $X$ can not lie on a Calabi-Yau.  Our
method is to describe explicitly {\em all} deformations of the
cone.

        We actually adopt the point of view of {\em extensions} of a
subvariety $V\subset\Pbb^n$ rather than deformations of the cone:
consider $W\subset\Pbb^{n+k}$ for which a codimension $k$
transversal linear section gives $V$.  Any projective variety is
an extension of its hyperplane section.  There is the notion of a
``universal extension'' of a given $V$ which in good cases is
governed by a Kodaira-Spencer map. Combining Theorem 2.8 with
Remark 1.10.1 gives

   \begin{theoA}  Suppose $V\subset\Pbb^n$ satisfies $(N_2)$,
       $V\not \cong \Pbb^1$, and
$W\subset\Pbb^{n+k}$ is an extension. From
the normal bundle sequence of $V$ in $W$ consider the coboundary map
   \[
\gamma : \Gamma (N_{V/W}(-1))\rightarrow H^1(\Theta_V(-1)).
   \]
Assume either
  \begin{enumerate}
    \item[(i)]   $\dim V =1$ and $\crk \gamma = \crk
\Phi(K_V,\Ocal_V(1))$ ($\Phi$ a Gaussian)
    \item[(ii)]  $\dim V \geq 2$ and $\gamma$ is an isomorphism.
  \end{enumerate}
Then $W$ is a {\bf universal} extension of $V$, in an appropriate
sense.
   \end{theoA}

As remarked in (1.10.4), this result may be used to deduce Mukai's
classification of Gorenstein Fano 3-folds of genus 7, 8, or 9 from
his earlier results on curves [M2].  In [BM], it is shown that for
a hyperplane section of a $K$-3 surface, $\gamma\neq 0$ (i.e., the
normal sequence is non-split) unless there is an involution fixing
the section.

Next we give a general method for constructing extensions of a
canonically polarized $(X,K_X)$, in case $X$ is a divisor on a
variety $Z$ (usually Fano).  We attribute this method to DuVal [D],
though a similar use of adjoint linear systems goes back to del
Pezzo.  We illustrate by considering a plane curve $C$ of degree
$d\geq $4 (cf. [E], [W2]).  Let $D$ be a smooth cubic
intersecting $C$ transversally, with $Z=C\cap D$.  The linear system
of curves of degree $d$ containing $Z$ defines a morphism on the
blow-up $B\rightarrow \Pbb^2$ of $Z$; this restricts to the canonical
embedding of the proper transform of $C$, and is an embedding except
for collapsing the proper transform of $D$ to a point.  The image of
$B$ is a normal surface with trivial dualizing sheaf and one simple
elliptic singularity, for which $C$ is a hyperplane section.  If
$d\geq 7$, keeping track of the normal sequence one can show all
extensions of $C\subset\Pbb^{g-1}$ arise from varying $D\in
|-K_{\Pbb^2}|$, with possibly worse (hence still non-rational)
singularities.  Since no such is a $K$-3, $C$ lies on no $K$-3.

This method generalizes greatly to large divisors on an
arbitrary $X$ (Corollary 3.12); one knows {\em all} extensions of the
canonical embedding.  Combining (3.12) and (3.13) yields the chief
result of this paper:

  \begin{theoB} Let $Z$ be a projective variety with isolated
Gorenstein singularities and $h^1(\Ocal_Z)=0$.  Then for sufficiently
ample smooth divisors $X$ on $Z$, $X$ does not sit on a Calabi-Yau,
unless $Z$ is already a Calabi-Yau.
  \end{theoB}

There is a much sharper statement when $\dim{Z}=2$ (Theorem 4.5).
We also remark that Theorem B is fairly easy in case $-K_Z$ has no
global sections; in that case the only extensions of $X$ are
cones.

  \medskip

\noindent {\bf Corollary 4.5:} {\em A smooth hypersurface of degree
$d$ in $\Pbb^n$ is a hyperplane section of a Calabi-Yau variety iff
$n+1\leq d\leq 2n+2$}.

  \medskip

In $\S$1, we introduce the (first) Kodaira-Spencer map of an
extension and prove (Theorem 1.9) that it can be used to prove a
given extension is universal.  We relate this map to the normal
bundle sequence of an extension in $\S$2.  In $\S$3, start with
$X\subset Y$ a divisor, with $h^1(\Ocal_Y)=0$; for an effective
divisor $D\subset Y$, we construct an extension of $X$ embedded
via $\Ocal_Y(X-D)\otimes \Ocal_X$, by using the rational map of
$Y$ given by sections of $\Gamma (Y,\Ical_{X\cap D}(X))$. One can
do this as well for a family of $D$'s; using $\S$2, one has
criteria that a universal extension has been constructed. These
results are applied to canonical embeddings of curves in Theorem
4.5, yielding the universal extension for all complete
intersection curves, most curves on the Hirzebruch surfaces
$\Fbb_n$, and bielliptic curves of genus $\geq 11$.  In
particular, we find that only the ``obvious''  complete
intersection curves actually lie on a $K$-3.  Finally, in $\S$5,
we give some necessary conditions for a smooth complete
intersection in $\Pbb^n$ of dimension at least 2 to be a
hyperplane section of a Calabi-Yau.

            It should be clear that a more accurate (though much more
cumbersome) title for this paper would have been ``Varieties which
are \emph{not} hyperplane sections of Calabi-Yau manifolds''.

        This research was partially supported by an NSF Grant.

\section{Extensions of subvarieties}

  \begin{numpar} %(1.1)
Let $V\subset\Pbb^n$ be a non-degenerate
subvariety.  Identify $\Pbb^n$ with a linear subspace $H$ of
$\Pbb^{n+k}$. A {\em ($k$-step) extension} of $V$ is a subvariety $W$
of $\Pbb^{n+k}$ so that $V=W\cap H$, transversally; thus, $\codim
(V,W)=k$, and the defining equations of $H$ form locally a regular
sequence in $W$ (cf. [Lv]).  Two extensions $(V,H,W,\Pbb^r)$ and
$(V,H',W',\Pbb^r)$ are {\em isomorphic} if there is an isomorphism of
$\Pbb^r$ sending $H$ onto $H'$, $W$ onto $W'$, and equal to the
identity on $V$; in particular, extensions $W$ and $W'$ of $V\subset
H\subset\Pbb^{n+k}$ are isomorphic if they differ by an automorphism
of $\Pbb^{n+k}$ which is the identity on $H$.  From
$W\subset\Pbb^{n+k}$, one may form a {\em sub-extension} $(V,H,W\cap
L,L)$, where $H\subset L\subset\Pbb^{n+k}$, as well as a {\em cone}
over $W$ in $\Pbb^{n+k+1}$.  When $V\subset\Pbb^n$ is linearly
normal, one may speak of an extension of $(V,\Ocal_V(1))$.
  \end{numpar}

  \begin{numpar}   %(1.2)
Consider $V\subset H\subset \Pbb^{n+1}$, and a projective cone $C(V)$
of $V$ over a point in $\Pbb^{n+1}-H$; so $C(V)\cap H=V$.  There is a
relation between extensions of $V\subset H$; part of the Hilbert
scheme of $C(V)\subset\Pbb^{n+1}$; and deformations of non-positive
weight of the affine cone of $V$ in $\Cbb^{n+1}$.  For instance,
suppose $W$ is a fibre in a deformation of $C(V)$ in $\Pbb^{n+1}$
which keeps fixed a hyperplane section: $W\cap H=C(V)\cap H=V$; then
$(V,H,W,\Pbb^{n+1})$ is an extension of $V$.  Conversely, if
$(V,H,W,\Pbb^{n+1})$ is a 1-step extension, then by a well-known
construction (e.g., [P]) $W$ degenerates in $\Pbb^{n+1}$ to a
projective cone over a hyperplane section $W\cap H=V$.  On the other
hand, it is a bit tricky to construct appropriate deformation
theories from the point of view of the Hilbert scheme with fixed
hyperplane section, or for the negative weight deformations of the
affine cone.
  \end{numpar}

  \begin{numpar}   %(1.3)
A $k$-step extension $W$ of $V$ is \textit{universal} if every
extension of $V$ is equivalent to a (possibly trivial) projective
cone over a unique
subextension. Such a $W$ has the weaker property (cf. [Lv]) that
any extension of it is a cone.  A universal extension need not
exist.  We give a criterion for a given extension of $V$ to be
universal, using a Kodaira-Spencer map.
  \end{numpar}

  \begin{numpar}   %(1.4)
  For $(V,H,W,\Pbb=\Pbb^{n+k})$ a $k$-step extension, there is a
natural {\em Kodaira-Spencer map}
  \[
KS: \Gamma (\Pbb , \Ical_H(1))^* \rightarrow
\Gamma (V,N_{V/H}(-1))/\Gamma (H,\Theta_H(-1)).
  \]
Here, $\Ical_H$ is the ideal sheaf of $H$ in $\Pbb$ and $N$
represents the normal sheaf. For, use the short
exact sequence
  \[
 0 \rightarrow \Gamma (\Ocal_H(1))^*
\rightarrow \Gamma (\Ocal_{\Pbb}(1))^*
\rightarrow \Gamma (\Ical_H(1))^* \rightarrow 0,
  \]
the sequence of maps
  \[
\Gamma (\Ocal_{\Pbb}(1))^* \cong \Gamma (\Theta_{\Pbb}(-1))
\rightarrow \Gamma (\Theta_{\Pbb |W}(-1))
\rightarrow \Gamma (N_{W/\Pbb}(-1)) \rightarrow \Gamma (N_{V/H}(-1)),
  \]
and the isomorphism
   \[
 \Gamma (\Ocal_H(1))^* \cong \Gamma (\Theta_H(-1)).
   \]
(The map $N_{W/\Pbb}\rightarrow N_{V/H}$  arises from
$\Ical_{W/\Pbb}\otimes\Ocal_H \cong \Ical_{V/H}$ because the
equations of $H$ give a regular sequence in $W$.)  The $KS$ map of a
sub-extension factors via the natural inclusion
   \[
\Gamma (L,\Ical_{H/L}(1))^* \subset \Gamma
(\Pbb ,\Ical_{H/\Pbb}(1))^*.
   \]
  \end{numpar}

One can define higher-order Kodaira-Spencer maps. For example, writing
$KS_{1}$ for the map above, one may consider
  \[
KS_{2}: \Ker KS_{1} \to \Gamma (N_{V/H}(-2)),
  \]
arising from considering the kernel of
  \[
N_{W/\Pbb}(-1) \to N_{V/H}(-1).
  \]

  \begin{examples}  %1.5
\item Suppose $Z\subset\Pbb^n = H$ is a smooth non-degenerate
subvariety, and $V=Z\cap Q$ is the transversal intersection with a
quadric hypersurface.  Choose coordinates $x_0,\cdots ,x_n$, $t$ on
$\Pbb^{n+1}$ so that $H=\{ t=0\}$.  Let $W\subset \Pbb^{n+1}$ be the
intersection of the projective cone over $Z$ with the quadric $\{
t^2=q(x_0,\cdots ,x_n)\}$, where $q$ defines $Q$.  Then $W$ is a
1-step extension of $V$; the Kodaira-Spencer map $KS$ is 0; $\Gamma
(N_{V/H}(-2))\neq 0$; and the higher-order map
  \[
\Ker KS \to \Gamma (N_{V/H}(-2))
  \]
is injective.
  \item Let $V=\{ F(x_0,\cdots ,x_n)=0\} \subset H=\Pbb^n$ be a
non-singular hypersurface of degree $d$.  Let $M_1,\cdots ,M_N$ be a
set of monomials giving a basis of the polynomials of degree $<d$ in
$\Cbb [x_0, \cdots ,x_n]/\Bigl( \frac{\partial F}{\partial
x_i}\Bigr)$.  Let $d_i=\deg M_i$.  Then in $\Pbb^{n+N}$, with
coordinates $x_0,\cdots ,x_n, t_1,\cdots ,t_N$, a universal extension
$V$ is given by
  \[
\Bigl\{ F+\sum t_i^{d-d_i} M_i =0\Bigr\} .
  \]
In this case, the Kodaira-Spencer maps are surjections, including
the $d-1$ higher-order ones, with target spaces $\Gamma
(N_{V/H}(-i))$, $i=2,\cdots ,d$.
  \end{examples}

  \begin{numpar}  %(1.6)
It will simplify matters to work with equations defining a
variety, so we will assume from now on that
  \begin{equation*}
\begin{matrix} V\subset H=\Pbb^n \mbox{ is normal and
projectively normal},\\
\mbox{with homogeneous coordinate ring } A=P/I(V)=P/I. \end{matrix}
   \tag{$\ast$}   \end{equation*}
This assumption guarantees that for an extension $W$ of $V$, elements
of $I(V)$ lift to $I(W)$.  Further, the target spaces of the Kodaira-Spencer maps are identified
from the following
well-known facts:
  \begin{enumerate}

  \item For all $i$, $\Gamma (V,N_{V/H}(i))=i^{\mbox{th}}$
graded piece of $\Hom_{P/I}(I/I^2, P/I)$.  \label{1.6.1}

  \item $\Gamma (\Theta_H(-1)) \to \Gamma (N_{V/H}(-1))$ is
given by differentiation.   \label{1.6.2}
  \end{enumerate}
  \end{numpar}

  \begin{numpar}  %(1.7)
Now suppose $V\subset H\subset \Pbb^{n+k}=\Pbb$, and choose coordinates
$x_0,x_1,\cdots, x_n,\break t_1, t_2,\cdots ,t_k$ so that $H$ is defined by
$\{ t_1=\cdots =t_k=0\}$.  By ($*$), a $k$-step extension
$W\subset\Pbb$ of $V$ gives a flat lifting of $I$ to $I(W)\subset
\Cbb [x_0,\cdots ,x_n,t_1,\cdots ,t_k]$.  The $KS$ map of $W$ maps
$\partial /\partial t_i$ in $\Gamma (\Ical_H(1))^*$ to a degree $-1$
homomorphism $I\to P/I$ \ref{1.6.1}: if $f\in I$, lift to an $F\in
I(W)$,
and map to the image $(\partial F/\partial t_i)\equiv g_i$ in $P/I$.
A coordinate change $x_{\gamma} \mapsto x_{\gamma} +\sum
a_{\gamma\delta} t_{\delta}$ ($\gamma =0,\cdots ,n; \delta =1,\cdots
,k$) replaces $g_i$ by $g_i+(\sum a_{\delta i}\partial /\partial
x_{\gamma})f$.  So $KS(\partial /\partial t_i)=[g_i]$ in $\Gamma
(N_V(-1))/\Gamma (\Theta_H(-1))$.  (Thus $KS$ is essentially the
Kodaira-Spencer map of a weight $-1$ deformation of the affine cone
$A=P/I$.)  Higher-order Kodaira-Spencer maps will involve higher-order
terms in variables $t_{i}$ in $\Ker{KS}$.
\end{numpar}

  \begin{numpar}  %(1.8)
  The following ought to be well known.
  \end{numpar}

  \begin{theorem} \label{th1.9} Let $V\subset\Pbb^n=H$ be a normal
and projectively normal subvariety, which satisfies
  \begin{equation*}
H^0(V, N_{V/H}(-2)) =0.  \tag{1.9.1}  \label{1.9.1}
  \end{equation*}
Then
  \begin{enumerate}
 \item[(a)] A $k$-step extension $W\subset\Pbb^{n+k}$ is uniquely
determined (up to equivalence) by its Kodaira-Spencer map.
 \item[(b)] A $k$-step extension $W$ with $KS$ an isomorphism is a
{\bf universal} extension, in that every extension is isomorphic to
a (cone over a)
unique sub-extension of $W$.
  \end{enumerate}
  \end{theorem}

  \begin{proof}  Choose coordinates as in (1.7).  Encode generators of
$I(V)$ in a $1\times \kappa$ row vector $f$.  Let $r$ be a
$\kappa\times\ell$ matrix of forms of $P$ whose columns generate the
module of relations for $f$.
  \begin{flalign*}
f\cdot r  &=0 \\
f\cdot r' &= 0, \mbox{ where $r'$ is $\kappa\times s$, implies
$r'=r\cdot b$, where $b$ is $\ell\times s$}.
  \end{flalign*}

Let $W\subset\Pbb^{n+k}$ be an extension of $V$.  $W$ is defined by
equations $F=f+\sum t^Ig^I$, where $I=(i_1,\cdots ,i_k)$,
$|I|=i_1+\cdots +i_k$, and $t^I = t_1^{i_1}t_2^{i_2} \cdots
t_k^{i_k}$;
relations are determined by $R=r+\sum t^Iu^I$, with $F\cdot R=0$.
Writing the $t$-linear part of $F$ as $\sum g_it_i$, the $KS$ map sends
$\partial /\partial t_i$ to the class of $g_i$ in $\Hom_{P/I}(I/I^2,
P/I)$ mod derivations.

Suppose $\bar{W}\subset \Pbb^{n+k}$ is another extension, with
equations and relations given by $\bar{F}, \bar{g}^I, \bar{R},
\bar{u}^I$, and with the same $KS$ map.  For each $i$ the classes of
$g_i$ and $\bar{g}_i$ are the same, so there are derivations
$D_i=\sum a_{i\alpha}\partial /\partial z_{\alpha}$ for which
  \[
g_i-\bar{g}_i = D_if + b_i\cdot f, \mbox{ where $b_i$ is a
$\kappa\times\kappa$ matrix}.
  \]
(Note that if the generators of $I$ have the same degree, then each
$b_i$ is 0.)  Replace $\bar{F}$ by $\bar{F}+\sum t_ib_i\cdot\bar{F}$,
then apply
the coordinate change $x_{\alpha}$ goes to $x_{\alpha}-\sum
a_{i\alpha}t_i$; one obtains an extension equivalent to $\bar{W}$,
with $g_i=\bar{g}_i$, all $i$.

Assume inductively that for all $I$ with $|I|\leq m$, $g^I=\bar{g}^I$
and for all $I$ with $|I| <m$, $u^I=\bar{u}^I$.  Writing the
coefficient of $g^I$ in the equations $F\cdot R=\bar{F}\cdot\bar{R}
=0$, one deduces
  \[
f\cdot u^I + (\cdots ) =0 \mbox{  and  } f\cdot\bar{u}^I + (\cdots
)=0,
  \]
where the expressions in parentheses are the same.  Therefore,
$f\cdot (u^I-\bar{u}^I)=0$, so that
  \[
u^I - \bar{u}^I = r\cdot b^I.
  \]
Thus, replacing $\bar{R}$ by $\bar{R}+\sum t^I\bar{R}\cdot b^I$ (sum
over all $I$ with $|I|=m$) gives another relation vector, now with
$u^I=\bar{u}^I$, all $|I|\leq m$.  Next fix an $I$ with $|I| =m+1$.
Again, one has equations
  \[
f\cdot u^I + (\cdots ) + g^I\cdot r = 0 \mbox{  and  }
f\cdot\bar{u}^I + (\cdots ) + \bar{g}^I\cdot r = 0,
  \]
where the expressions in parentheses are equal to each other, by the
inductive hypothesis.  Thus,
  \[
(g^I-\bar{g}^I)\cdot r = f\cdot v^I,
  \]
whence $g^I-\bar{g}^I$ defines a homomorphism of $I$ into $P/I$ of
degree $-(m+1)$.  But the hypothesis $\Gamma (N_{V/H}(-2))= 0$
implies also that $\Gamma (N_{V/H}(-i))=0$, $i\geq 2$, whence (1.6.1)
a homomorphism $I\to P/I$ of degree $\leq -2$ must be 0. Thus, one
can write
  \[
g^I-\bar{g}^I = f\cdot w^I.
  \]
Replacing $\bar{F}$ by $\bar{F}+\sum t^I\bar{F}\cdot w^I$ (over all
$I$ with $|I|=m+1$), one may assume $g^I=\bar{g}^I$, all such $I$, so
the induction may continue.  In this way, since the sum is finite by
degree considerations, one may conclude that $W$ and $\bar{W}$ are
equivalent.

Next, fix a $k$-step extension $W$ with $KS$ an isomorphism;
represent it as above by equations
  \[
F = f+\sum t_ig_i + \sum t^Ig^I.
  \]
Consider an arbitrary $j$-step extension $W'$, with linear section
defined by the vanishing of coordinates $s_1, s_2, \cdots ,s_j$.
Since the classes $[g_i]$ form a basis, one can after allowable
linear coordinate changes in $\Pbb^{n+j}$ write the equations for
$W'$ as
  \[
\bar{F} = f+\sum' s_ig_i + \sum s^Ih^I,
  \]
where the first sum involves only the first $j'$ variables (some
$j'\leq j$).  Using $W$, one may write down another extension in
$\Pbb^{n+j}$ with the same $KS$ map; it is the cone over the
subextension in $\Pbb^{n+j'}$
  \[
F^* = f + \sum' s_ig_i + \sum' s^Ig^I,
  \]
where both sums involve only the first $j'$ variables.  By (a), this
extension is isomorphic to $W'$.
  \end{proof}

  \smallskip

  \begin{remarks}
   \item \label{1.10.1} It is proved in ([W1], (2.8)) that the
vanishing condition (1.9.1) is automatic when $V\subset \Pbb^n$ is
smooth and satisfies ($N_2$):  that is, when $I(V)$ is generated by
quadrics and the relations are generated by linear ones.

  \item \label{1.10.2} Presumably one could prove a universality
result for an extension $V\subset W$ without assuming (1.9.1).  One
should require then that if $\Gamma(N_V(-i))=0$, all $i>d$, then the
Kodaira-Spencer maps $KS_{i}$ are surjective for all $i<d$, and an
isomorphism for $i\geq d$.  (Compare (1.5.2).)  However, this would need
to be worked out.

\item \label{1.10.3} For $n\geq3$ the Segre embedding $\Pbb^{1}\times\Pbb^{n-1}\subset\Pbb^{2n-1}$ gives an extension of the twisted rational curve of degree $n$; it is universal for $n\neq4$, but the Kodaira-Spencer map is not surjective once $n>3.$

  \item \label{1.10.4} This Theorem allows a somewhat simpler
proof of $S$.  Mukai's classification [M2] of Gorenstein Fano 3-folds
$X$, with $-K$ a generator of Pic($S$), and with $g=7,8,$ or 9: they
are linear sections of an appropriate homogeneous space.  Mukai shows
that for the anticanonical embedding of $X$, the general linear
section of codimension 2 is a canonical curve $C$, which (by his
earlier work)  is a linear section of a certain $G/P\subset\Pbb^N$.
One claims that $X$ is also a linear section of $G/P$, as will follow
from 1.9 applied to the extension $C\subset G/P$.  That the canonical
embedding of $C$ satifies (1.9.1) follows, e.g., from (4.1.2) below
plus results of [M1].  $G/P$ is smooth, hence not a cone, so the $KS$
map of this extension of $C$ is injective (by 1.9).  But $KS$ is in
fact bijective, by a dimension count; $\dim (G/P)-1$ is the corank of
the Gaussian map --- cf. [CM], ([W2], (6.5)).

  \end{remarks}

  \section{$KS$ and the normal bundle sequence}

  \begin{numpar}  We frequently compose the
  Kodaira-Spencer map with a coboundary, because of the simple
  \end{numpar}
  \begin{lemma}  Let $V\subset H=\Pbb^n$ be a smooth linearly
  normal subvariety, and $V\neq\Pbb^{1}$.  Then the coboundary map associated to the normal
  bundle sequence
  \[
  \delta: \Gamma (V,N_{V/H}(-1))/\Gamma (H,\Theta_H(-1)) \rightarrow
H^1(\Theta_V(-1))
  \]
is injective.

 \end{lemma}
 \begin{proof}  By the usual long exact sequence in cohomology of the
 normal bundle sequence for $V\subset H$, it suffices to show that
 \[
 \Gamma(\Theta_H(-1))\cong\Gamma(\Theta_H(-1) \otimes \Ocal_V).
 \]
By the standard presentation of $\Theta_H$, this isomorphism is
equivalent to the injectivity of a map
\[
H^1( \Ocal_V(-1)) \rightarrow H^1(\Ocal_V^{n+1}).
\]
When $\dim(V)\geq 2$, the first group is $0$ by Kodaira vanishing.  When
$\dim(V)=1$, by duality (and linear normality) one needs the
surjectivity on the curve $V$ of
  \[
\Gamma(K_V) \otimes \Gamma(\Ocal(1)) \rightarrow \Gamma(K_V(1)),
  \]
which is a well-known result of Petri when $V\neq\Pbb^{1}$.
\end{proof}
\begin{remarks}  \item Similarly, one can show that the coboundary $\delta$ is an
    isomorphism when $H^1(\Ocal_V)=0$ and either $\dim(V)\geq 3$, or $\dim(V)=2$
    and $\Gamma(K_V) \otimes \Gamma(\Ocal(1)) \rightarrow \Gamma(K_V(1))$ is
    surjective.
    \item We conclude that the Kodaira-Spencer map for an extension of a
    linearly normal subvariety $V\subset \Pbb^n$ has rank at most
    $h^1(\Theta_V(-1))$.
\end{remarks}

 \begin{numpar}  We can now relate the Kodaira-Spencer map of an extension
to the class of the normal bundle sequence of $V\subset W$.  Assume $V\subset H=\Pbb$ is smooth and projectively normal, and
$(V,H,W,\Pbb^{n+k})$ is a k-step extension.

  \end{numpar}

  \begin{numpar} The
normal bundle of $V$ in $W$ is the restriction of the normal bundle
of $H\subset\Pbb$; tensoring the surjection
  \[
\Gamma (\Ical_H(1)) \otimes \Ocal_{\Pbb}(-1) \twoheadrightarrow
\Ical_H
  \]
with $\Ocal_H$, dualizing, and restricting, one has
  \[
N_{V/W} \cong \Gamma (\Ical_H(1))^* \otimes \Ocal_V(1).
  \]
So, one has an identification
\[\Gamma(N_{V/W}(-1))\cong\Gamma(\Ical_H(1))^*.\]
The normal bundle sequence of $V$ in $W$
  \[
0 \to \Theta_V \to \Theta_W\otimes\Ocal_V \to N_{V/W} \to 0
  \]
is therefore determined as a bundle extension by a map
  \[
\gamma: \Gamma (\Ical_H(1))^* \to H^1(\Theta_V(-1)).
  \]
  \end{numpar}

  \begin{proposition} Let $V\subset \Pbb^n=H$ be smooth and
projectively normal ($V\not\cong\Pbb^1$), $W\subset\Pbb^{n+k} = \Pbb$ a
$k$-step extension.  Then the bundle extension map $\gamma$ above is
the negative of the composed map $\delta \cdot KS$.
  \end{proposition}

  \begin{proof}  Choose coordinates as before, with $I\subset P$,
and let $\bar{I} = I+(t_1,\cdots ,t_k)\subset \bar{P} = P[t_i]$.
From the exact sequence
  \begin{enumerate}
\item \label{1.12.1} $\qquad\qquad 0 \to N_{V/H} \to N_{V/\Pbb} \to
N_{H/\Pbb} \otimes \Ocal_V \to 0$
  \end{enumerate}
and \ref{1.6.1}, one deduces that $\Gamma (N_{V/\Pbb}(-1))$ consists
of degree $-1$ homomorphisms $\bar{I}/\bar{I}^2 \to \bar{P}/\bar{I} =
P/I$.  The subspace $\Gamma (N_{V/H}(-1))$ consists of homomorphisms
for which the $t_i$ go to 0.

For an extension $W\subset\Pbb$, each $f\in I$ lifts to $F\in
I(W)\subset \bar{I}$, where $F=f+\sum t_ig_i + \cdots$.  Thus,
$KS(\partial /\partial t_{\mu})$ is the class of the homomorphism in
$\Gamma (N_{V/\Pbb}(-1))$ sending $f$ to $g_{\mu}$, and every $t_i$
to 0.

The map $\gamma$ factors through $\Gamma (\Ical_H(1))^*\cong
\Gamma (N_{V/W}(-1))\subset \Gamma (N_{V/\Pbb}(-1))$ (an extension
provides a splitting of {2.6.1}).  Now, $\partial /\partial
t_{\mu}$ sends $t_{\nu}$ to $\delta_{\nu\mu}$; since
  \[
f = -\sum t_ig_i - \cdots
  \]
on $W$, $\partial /\partial t_{\mu}$ sends $f$ to $-g_{\mu}$.

Thus $(\gamma + \delta \cdot KS)(\partial /\partial t_{\mu})$ is
represented in $\Gamma(N_{V/\Pbb}(-1))$ by a map sending each $f\in I$ to 0 and $t_i$ to
$\delta_{i\mu}$.  This is exactly the image of the element
``$\partial /\partial t_{\mu}$'' via
  \[
\Gamma (\Theta_{\Pbb}(-1)) \to \Gamma
(\Theta_{\Pbb}\otimes\Ocal_V(-1)) \to \Gamma (N_{V/\Pbb}(-1)).
  \]
Since the cokernel of the last map is contained in
$H^1(\Theta_V(-1))$, we have the assertion of the Proposition.
  \end{proof}

  \begin{remark} Let us denote the composed Kodaira-Spencer map
$\delta \cdot KS$ by $\overline{KS}$; this map has target space
$H^1(\Theta_V(-1))$.  We shall also abuse notation and not
distinguish between $\overline{KS}$ and its negative $\gamma$.
  \end{remark}
\begin{theorem} Suppose $V\subset\Pbb^n$ is a projectively normal
embedding of a smooth variety ($V\not\cong\Pbb^{1}$), for
which $\Gamma(N_{V/\Pbb^n}(-2))=0$.  Let
$W\subset \Pbb^{n+k}$ be an extension. From
the normal sequence of $V$ in $W$ consider the coboundary map
   \[
\gamma : \Gamma (N_{V/W}(-1))\rightarrow H^1(\Theta_V(-1)).
   \]
Assume either
  \begin{enumerate}
    \item[(i)]   $\dim V =1$ and $\crk \gamma = \crk
\Phi(K_V,\Ocal_V(1))$ ($\Phi$ a Gaussian)
    \item[(ii)]  $\dim V \geq 2$ and $\gamma$ is an isomorphism.
  \end{enumerate}
Then $W$ is a {\bf universal} extension of $V$.
\end{theorem}
 \begin{proof} By Theorem 1.9, one must show that $KS$ is an isomorphism;
composing with the injective coboundary map $\delta$, via (2.6) it
  suffices to consider $\gamma$.  It remains necessary only to
     add that for a projectively normal curve, the coboundary map $\delta$
     has the same corank as the Gaussian $\Phi(K_V,\Ocal_V(1))$, by e.g.
     [CM], (1.2).
\end{proof}

  \section{A Du Val-type construction of extensions}

  \begin{numpar} For $E\subset X$ a smooth Cartier divisor on a
projective variety, the {\em normal class} is the element in
$H^1(\Theta_E(-E))$ given by extension class of the normal bundle
sequence
  \[
0 \to \Theta_E \to \Theta_X\otimes\Ocal_E \to \Ocal_E(E) \to 0.
  \]
Let $Z\subset E$ be an effective Cartier divisor, $\pi : \tilde{X} =
Bl_ZX \to X$ be the blow-up of the ideal sheaf $\Ical_Z$ of $Z$ in
$X$, $\tilde{Z}\subset\tilde{X}$ the exceptional divisor, and
$\tilde{E}\subset\tilde{X}$ the proper transform of $E$.  Via the
isomorphism
  \[
\pi |_{\tilde{E}}: \tilde{E} \cong E,
  \]
one can identify the normal bundles
  \[
N_{\tilde{E}/\tilde{X}} \cong N_{E/X}(-Z).
  \]
In particular, the normal sequences are comparable:
  \begin{alignat*}{3}
  0 \to &\Theta_{\tilde{E}} \to &\Theta_{\tilde{X}} &\otimes
\Ocal_{\tilde{E}} \to &\Ocal_E &(E-Z) \to 0 \\
 &\| &&\cap &\;\cap \\
0 \to &\Theta_E \to &\Theta_X &\otimes\Ocal_E \to &\Ocal_E &(E)
\rightarrow 0.
  \end{alignat*}
Thus, the normal classes of $E\subset X$ and $\tilde{E}\subset
\tilde{X}$ are related via
  \[
H^1(\Theta_E(-E)) \to H^1(\Theta_E(-(E-Z))) =
H^1(\Theta_{\tilde{E}}(-\tilde{E})).
  \]
  \end{numpar}

  \begin{numpar}
For a line bundle $L$ on $X$ one has an exact sequence
  \[
0 \to L(-E) \to \Ical_Z\cdot L \to L\otimes\Ocal_E(-Z) \to 0.
  \]
If the linear system associated to $W\equiv \Gamma (\Ical_ZL)\subset
\Gamma (L)$ has base scheme exactly $Z$, then the rational map $X- \,
- \, \to \Pbb (W^*)$ becomes a morphism
  \[
\rho : \tilde{X} \to \Pbb (W^*)
  \]
associated to the line bundle $\pi^*(L)(-\tilde{Z})$.  The linear
system defining $\rho|_{\tilde{E}}$ is given by the image of $W$ in
$\Gamma (E,L\otimes\Ocal_E(-Z))$.
  \end{numpar}

  \begin{numpar} Assuming $L=\Ocal (E)$, $\Ocal_E(E-Z)$ is very
ample, and $H^1(\Ocal_X)=0$, one has the exact sequence
  \begin{enumerate}
\item \label{3.3.1} $\qquad\qquad 0\to \Gamma (\Ocal_X) \to \Gamma
(\Ical_Z\Ocal (E)) \to \Gamma (E, \Ocal_E(E-Z)) \to 0.$
  \end{enumerate}
It follows that $Z$ is the base scheme associated to $\Gamma
(\Ical_ZL)$.  Note $H^1(\Ocal_{\tilde{X}})=0$ and
$\Ocal_{\tilde{E}}(\tilde{E}) \cong \Ocal_E(E-Z)$ is very ample.
  \end{numpar}

  \begin{lemma} \label{lem3.4} Let $A$ be a smooth divisor on a
projective variety $B$ with $H^1(\Ocal_B)=0$.  Suppose $\Ocal_A(A)$
is very ample.  Then $\Ocal (A)$ is basepoint-free and gives a
birational morphism $\rho : B\to\Pbb (\Gamma (B,\Ocal (A))^*)$.
$\rho$ is an isomorphism in a neighborhood of $A$, and maps $A$
isomorphically to a hyperplane section $\bar{A}$ of $\rho (B) \equiv
\bar{B}$.
  \end{lemma}

  \begin{proof}  The morphism associated to $\Ocal (A)$ clearly
separates points of $A$ from points off $A$ and, as $\Gamma (\Ocal
(A))\twoheadrightarrow \Gamma (\Ocal_A(A))$, from other points on
$A$.  By hypothesis tangent vectors on $A$ are separated; since $A$
is smooth, tangent vectors of $X$ at points of $A$ are separated.
  \end{proof}

  \begin{proposition} \label{prop3.5} Let $E\subset X$ be a smooth
Cartier divisor on a projective variety with $h^1(\Ocal_X)=0$,
$Z\subset E$ a divisor with $\Ocal_E(E-Z)$ very ample.  Then
  \begin{enumerate}
   \item[(i)] the linear system $\Gamma (\Ical_Z\Ocal (E))$ gives a
(1-step) extension $\bar{X}$ of $(E, \Ocal_E({E-Z}))$
    \item[(ii)] the normal class of $E\subset\bar{X}$ is the image of
the normal class of $E\subset X$ via $H^1(E,\Theta_E(-E))\to
H^1(E,\Theta_E(-(E-Z)))$, or arising from
  \begin{align*}
 &\Ocal_E \\[-2ex]
 &\,\cap \\[-2ex]
0 \to \Theta_E(-(E-Z)) \to \Theta_X\otimes\Ocal_E(-(E-Z)) \to
\Ocal_E &(Z) \to 0.
  \end{align*}
  \end{enumerate}
  \end{proposition}

  \begin{remarks}  \item In case $X=\Pbb^2$, $E$ is a smooth cubic,
and $Z\subset E$ is a scheme of length $\leq 6$, the construction
above gives a Del Pezzo surface, viewed as an extension of $E$.
  \item The sequence \ref{3.3.1} shows that the domain of the $KS$
map for the extension in (3.5.ii) is isomorphic to $\Gamma
(\Ocal_X)^*$, hence has a natural basis element.  Thus, in case
$\Ocal_E(E-Z)$ is projectively normal, the $\overline{KS}$ map is
described by the (negative of the) normal class of
$E\subset\bar{X}$.
  \end{remarks}

  \begin{corollary} \label{cor3.7} Let $E$ and $F$ be effective Cartier
divisors on a projective variety $X$, with $E$ smooth and
$h^1(\Ocal_X)=0$.  Suppose $\Ocal_X(E-F)$ is very ample.  Then
denoting $Z=E\cap F$:
  \begin{enumerate}
 \item \label{3.7.1} The linear system of $\Gamma (\Ical_Z\Ocal (E))$
has base scheme $Z$ and defines an isomorphism off $\Supp (F)$.

  \item \label{3.7.2} $\rho :\tilde{X} = Bl_ZX\to \Pbb (\Gamma
(\Ical_Z\Ocal (E))^*)$ is an isomorphism off $\tilde{F}$, the proper
transform of $F$, and collapses $\tilde{F}$ to a 0-dimensional
scheme.

  \item $E\subset\rho (\tilde{X}) =\bar{X}\subset\Pbb (\Gamma
(\Ical_Z\Ocal (E))^*)$ is an extension of $E$ embedded by the very
ample line bundle $\Ocal_E(E-F)$.  \label{3.7.3}

  \item The $\overline{KS}$ map (2.7) for the extension $E\subset\rho (\tilde{X})$ is
given by the composition  \label{3.7.4}
  \[
 \Gamma (\Ocal_E) \to \Gamma (\Ocal_E(F)) \to
H^1(E,\Theta_E(-(E-F))).
  \]
  \end{enumerate}
  \end{corollary}

  \begin{proof}  Since the linear system of $\Gamma (\Ical_Z \Ocal
(E))$ is generated by $E$ and $F+D$, where $D$ runs through the very
ample system $|E-F|$, \ref{3.7.1} follows.

For \ref{3.7.2}, let $U=\mbox{Spec }R\subset X$ be an affine open
neighborhood, with $E, F$ defined locally by $f,g\in R$,
respectively.  On $U$, the rational map may be written as
  \[
[f,g, gh_1, \cdots ,gh_s],
  \]
where $(h_1, \cdots ,h_s): U\to \Cbb^s$ is an embedding.  The
complement of $\tilde{F}$ in $\pi^{-1}(U)$ is an affine
$\tilde{U}\subset\tilde{X}$, with coordinate ring $R[t]/f-tg$, with
$\tilde{E}$ defined by $t, \tilde{Z}$ by $g$.  On $\tilde{U}$, the
rational map becomes a morphism
  \[
[t,1,h_1,\cdots ,h_s];
  \]
as this gives an embedding on $\Cbb\times U$, which contains
$\tilde{U}$, the first assertion follows.  Note also that
$\Ical_Z\Ocal (E)$ restricted to $F$ is $\cong \Ocal_F$.  The
remaining claims follow from Proposition 3.5.
  \end{proof}

  \begin{numpar}  One may form $k$-step extensions in the DuVal-type
construction above by considering the complete linear system of $Z$
on $E$.  Let $Y$ be a projective variety with $h^1(\Ocal_Y)=0$,
$D\subset Y$ a smooth Cartier divisor, $M$ a line bundle on $D$ so
that $\Ocal_D(D-M)$ is very ample, and $V\subset\Gamma (D,M)$ a non-0
linear subspace.  The linear system associated to $V$ is represented
by $Z\subset D\times\Pbb (V)$, a relatively effective Cartier divisor
over $\Pbb (V)$.  In the previous notation, let
  \end{numpar}
  \begin{align*}
X &= Y\times \Pbb (V) \\
E &= D\times \Pbb (V)   \\
Z &\subset D\times \Pbb (V) \mbox{ as above }\\
L &= \Ocal_Y(D)\boxtimes \Ocal (1).
  \end{align*}
Note that $\Ocal_E(Z) \cong \Ocal_D(M)\boxtimes\Ocal (1)$.  The
restriction of $\Ical_ZL$ to $E$ is
  \[
L\otimes\Ocal_E(-Z) \cong \pi_1^*\Ocal_D(D-M).
  \]
Since $H^1(X,L(-E)) \cong H^1(Y,\Ocal_Y)\otimes H^0(\Pbb (V), \Ocal
(1))=0$, one has the exact sequence
  \[
0 \to \Gamma (\pi^*_2\Ocal (1)) \to \Gamma (\Ical_ZL) \to \Gamma
(\pi_1^*\Ocal_D(D-M)) \to 0.
  \]
The first subspace is canonically $V^*$, and yields divisors of the
form $E+\pi_2^{-1}H$, where $H\subset\Pbb (V)$ is a hyperplane; so
$\Gamma (\Ical_ZL)$ has no base points off $E$.  $\Ocal_D(D-M)$ is
very ample, hence free, so there are no base points on $E$ except
along $Z$.  The map $\rho : \tilde{X} \to \Pbb (\Gamma
(\Ical_ZL)^*)=\Pbb^N$ restricted to $\tilde{E}\cong E=D\times\Pbb
(V)$ is the projection map to $D$ followed by the embedding given by
$\Ocal_D(D-M)$.  Intersecting the linear space $\Pbb
(V)\subset\Pbb^N$ with $\rho (\tilde{X})$ and pulling back to
$\tilde{X}$ is (as a scheme) the intersection of the divisors
$E+\pi_2^{-1}H$, which equals $E$.  We conclude that the smooth
variety $\rho (\tilde{E})$ equals $\rho (\tilde{X})\cap\Pbb (V)$.  In
particular, $\rho$ is birational onto its image, which is a $(\dim
V)$-step extension of $\rho (\tilde{E})\cong D$, embedded via
$\Ocal_D(D-M)$.  Note finally that $V$ is naturally isomorphic to
$\Gamma (\Ical_{\Pbb (V)}(1))^*$.  We change notation and summarize
in the

  \begin{theorem}
Let $E\subset X$ be a smooth Cartier divisor on a projective variety
with $h^1(\Ocal_X)=0$, $M$ a line bundle on $E$ with $\Ocal_E(E-M)$
very ample.
  \begin{enumerate}
    \item[(a)] Then for every non-0 $\;V\subset \Gamma (E,M)$, there
is
a $(\dim V)$-step extension $\bar{X}$ of $E$ embedded by
$\Ocal_E(E-M)$.
    \item[(b)] Assume further $\Ocal_E(E-M)$ is projectively normal.
Then the $\overline{KS}$ map of the extension $\bar{X}$ is the
composition
  \[
V\subset\Gamma (E,M) \to H^1(E,\Theta_E(-E)\otimes \Ocal (M)),
  \]
the second map arising from the normal sequence of $E\subset X$.
   \end{enumerate}
   \end{theorem}

\begin{corollary} \label{cor3.10}
 Let $E\subset X$ be a smooth Cartier divisor on a
normal projective variety with $h^1(\Ocal_X)=0$.  Assume $M$ a line
bundle on $X$ so that
  \begin{enumerate}
 \item[(a)] \quad $\Ocal_X(E)\otimes M^{-1}$ is very ample
 \item[(b)] \quad $\Ocal_X(E)\otimes M^{-1}\otimes \Ocal_E$ is projectively normal.
  \end{enumerate}
Then there is a $h^0(M)$-step extension of $(E,\Ocal_E(E)\otimes
M^{-1})$, whose 1-step sub-extensions are formed by blowing up $X$ at
$E\cap F$ (where $F$ is effective, $\Ocal (F) \cong M$), and blowing
down each component of $\tilde{F}$.

If further
  \[
H^0(\Theta_X\otimes\Ocal_E(-E)\otimes M))=0,
  \]
then the Kodaira-Spencer map of the extension is injective.
  \end{corollary}

  \begin{proof} By the Kodaira vanishing theorem for a normal
variety, $h^i(X,\Ocal_E(-E)\otimes M)=0$, $i=0,1$, so $\Gamma
(X,M)\cong \Gamma (E,M\otimes \Ocal_E)$.  The normal sequence of $E$
in $X$ gives the exact
  \[
H^0(E,\Theta_X\otimes\Ocal_E(-E)\otimes M) \to H^0(E,M\otimes\Ocal_E)
\to H^1(E,\Theta_E(-E)\otimes M).
  \]
Now all assertions follow from the Theorem.
  \end{proof}

  \begin{remarks}
   \item An extension constructed as in Corollary \ref{cor3.10}
has a non-rational singularity if $h^{n-1}(\Ocal_F)\neq 0$ ($n=\dim
X$).  For, one blows up $X$ along $F\cap E$ and blows down $\tilde{F}
\cong F$ via $\rho$, to a finite set of points.  In the notation
above, $R^{n-1}\rho_*\Ocal_{\tilde{X}} \twoheadrightarrow
H^{n-1}(\Ocal_{\tilde{F}}) \neq 0$.

    \item Suppose $(E,\Ocal_E(E)\otimes M^{-1})$ satisfies ($N_2$)
(or just (\ref{1.9.1})).  A dimension count plus Corollary 3.10 can be used
to show the constructed extension is
universal.  This is of special interest when $M=-K_X$.
  \end{remarks}

  \begin{corollary} \label{cor3.12} Let $E\subset X$ be a smooth
Cartier divisor on a normal Gorenstein projective variety with
$h^1(\Ocal_X)=0$.  Assume
  \begin{enumerate}
   \item[(a)] \quad $K_X+E$ is very ample, and $K_E$ is projectively
normal
   \item[(b)] \quad $h^0(E,\Theta_X(-(K_X+E))\otimes \Ocal_E)=0$
   \item[(c)] \quad For the embedding $E\subset\Pbb = \Pbb (\Gamma
(K_E)^*)$ given via $K_E$,
    \begin{align*}
\dim \Gamma &(N_{E/\Pbb}(-1))/ \Gamma (\Theta_{\Pbb}(-1)) = h^0(-K_X)
\\
  \Gamma &(N_{E/\Pbb}(-2)) = 0.
     \end{align*}
  \end{enumerate}
Then the construction in (3.8) gives a universal extension of
($E,K_E$).  Any non-conical 1-step extension of $(E,K_E)$ is constructed
from an anti-canonical divisor $F$ on $X$, by blowing up $E\cap F$,
and blowing down $\tilde{F}$.  In particular, if $K_X\not\cong\Ocal_X$,
then every extension of $(E,K_E)$ has a non-rational
singularity.
  \end{corollary}

  \begin{proof}  If $h^0(-K_X)=0$, by (c) the $KS$ map of any
extension is 0; by Theorem \ref{th1.9}, it must be a cone
over $(E, K_E)$ (necessarily non-rational, as $K_E$ ample implies
$h^{n-1}(\Ocal_E)\neq 0$).

If $h^0(-K_X)\neq 0$, Corollary \ref{cor3.10} applied to $M=-K_X$ gives
an extension of $(E,K_E)$ with $KS$ an isomorphism (via (c)).
By (\ref{th1.9}), it is universal, and every extension is obtained by
blowing up $X$ along
$E\cap F$, where $F\in |-K_X|$ is an anticanonical divisor.  If
$K_X\cong \Ocal_X$, one obtains that $X$ is the essentially unique
extension of $(E,K_E)$.  If $K_X\cong\Ocal_X(-F)$, with $F\neq 0$,
then
  \[
0 \to K_X \to \Ocal_X \to \Ocal_F \to 0
  \]
and the equalities $h^n(\Ocal_X)=h^0(K_X)=0$, $h^n(K_X)=1$ imply
that $h^{n-1}(\Ocal_F)\neq 0$.  Applying Remark (3.11.1), one has a
non-rational singularity on every 1-step extension.
  \end{proof}

  \begin{remark}  It is important to note that the conditions of
  Corollary 3.12 hold for sufficiently ample $E$ on $X$; that is,
for $D$ is ample and $n$ sufficiently large, a smooth $E$ in
$|nD|$ satisfies (a)--(c).  This is clear for (a) and (b),
since $h^1(\Ocal_X)=0$ implies the surjectivity of
\[
\Gamma (X,K_X+E) \rightarrow \Gamma (E,K_E).
\]
For (c), it follows e.g. from [I] that for large $E$, $K_X+E$
satisfies the syzygy condition ($N_2$).  So it remains to check
the dimension claim.  If $\dim{X}\geq 3$, by (2.3.2) we need only
show that for large E
\[
\dim{H^1(\Theta_E(-K_E))}=h^0(-K_X).
\]
The normal bundle sequence for $E\subset X$ makes clear that for
$E$ large the first term is $h^0(E,-K_{E})$, which in turn equals
the second term for $E$ large.  When $\dim{X}=2$, the result is
well-known (e.g., [DM] or (4.4) below).

  \end{remark}

  \section{The case of curves}

  \begin{numpar}  If $C$ is a smooth non-hyperelliptic curve of genus
$g\geq 3$, then $K_C$ is very ample, and the canonical embedding
$C\subset\Pbb^{g-1}=\Pbb$ is projectively normal (Noether's theorem).
An extension $W\subset\Pbb^g$ of $C$ is a normal Gorenstein surface
with $K_W\cong\Ocal_W$ and $h^1(\Ocal_W)=0$ (a ``canonically trivial
surface.'')  To apply Corollary \ref{cor3.12}, note
  \begin{enumerate}
   \item \label{3.1.1} $\dim H^0(N_{\Pbb}(-1))/H^0(\Theta_{\Pbb}(-1))
= \crk\; \Phi_K: \Lambda^2\Gamma (K) \to \Gamma (K^{\otimes 3})  =
$ corank of Gaussian-Wahl map [W1].
   \item \label{3.1.2} $H^0(N_{C/\Pbb}(-i))=0$, $i\geq 2$, unless $C$
is trigonal of genus $\leq 10$; the intersection of a del Pezzo
surface in $\Pbb^{g-1}$ with a quadric ($g\leq 10$); or is
bielliptic.
  \end{enumerate}
(4.1.2) asserts the Gaussian $\Phi(K,K^{\otimes i})$ is surjective
for $i\geq 2$ except in these cases.  When the Clifford index $\geq
3$, one has property $(N_2)$ (by [S], [V]); so recall \ref{1.10.1}.
The other cases can be deduced via [CM]; [T], p. 161; and [W2],
$\S$5.8.
  \end{numpar}

  \begin{numpar}  Conditions (b) and (c) in Corollary \ref{cor3.12}
arise dually in the
  \end{numpar}

  \begin{proposition} \label{prop3.3} (cf. [DM]).  Let $C\subset X$
be a smooth Cartier divisor on a normal Gorenstein surface with
$h^1(\Ocal_X)=0$.  Assume
  \begin{enumerate}
 \item[(a)] \quad $H^1(\Omega^1_X (2K_X+2C)\otimes\Ocal_C) = 0$
 \item[(b)] \quad $K_X+C$ is ample
 \item[(c)] \quad The Gaussian $\Phi_{K_X +C}$ is surjective
 \item[(d)] \quad $H^1(\Omega^1_X (2K_X +C))=0$.
  \end{enumerate}
Then
  \[
  \corank \Phi_K = h^0(-K_X),
  \]
i.e., ``$X$ computes the Gaussian of $C$.''
  \end{proposition}

  \begin{proof}  A straightforward diagram chase as in [DM], Lemma
2.6, except that one uses $h^i(X,\Ocal_X(-(K_X+C)))=0$ ($i=0,1$), by
Kodaira Vanishing for a normal surface.
  \end{proof}

  \begin{numpar} Of course sufficiently ample $C$ on $X$ satisfy
conditions (a)--(d) of the Proposition.  We paraphrase Corollary
\ref{cor3.12} as follows:
  \end{numpar}

 \begin{theorem} \label{th3.5}
Suppose $X$ is a Gorenstein projective surface with $h^1(\Ocal_X)=0$.
Let $C\subset X$ be a smooth Cartier divisor, not one of the special
curves of genus $\leq 10$ of \ref{3.1.2}, nor bielliptic.  Assume
   \begin{enumerate}
 \item[(i)] $K_X+C$ is very ample
 \item[(ii)] $X$ computes the Gaussian of $C$.
   \end{enumerate}
Then every extension of the canonical embedding of $C$ is constructed
as in $\S$3 from an anti-canonical divisor on $X$.  In particular,
$C$ does not sit on a $K$-3 surface, unless $X$ is a $K$-3 surface.
  \end{theorem}

  \begin{remark} In the moduli space $\Mcal_g$, there are generally
many irreducible components of the locus of curves with Gaussian
corank 1, besides the $(g+19)$-dimensional stratum of $K$-3 curves.
For instance, let $X\subset\Pbb^{g'}$ be the cone over a general
canonical curve $C'$ of genus $g'\geq 12$.  By the Theorem, a smooth hypersurface
section $C$ of degree $d$ large has genus $g=d^2(g' -1)+1$ and
Gaussian corank 1, and uniquely determines $X$ (hence $C$).  Varying
$C'$ and the hypersurface section gives a $(g+2g'-3)$-dimensional
locus in $\Mcal_g$, with Gaussian corank 1.
  \end{remark}

  \begin{numpar} One can say precisely which complete intersection
curves sit on $K$-3 surfaces, and more generally describe the
canonically trivial surfaces of which they are hyperplane sections.
Let $C\subset\Pbb^n$ be a complete intersection curve of multidegree
$2\leq d_1\leq d_2 \cdots \leq d_{n-1} \equiv d$; assume $g\geq 2$,
so $\kappa =\sum d_i-(n+1)>0$ and $K_C=\Ocal_C(\kappa )$.  Let $X$ be
a general complete intersection surface of multidegree $d_1\leq
d_2\leq \cdots \leq d_{n-2}$ which contains $C$ ($X$ is unique unless
$d_{n-2}=d_{n-1}$); then $K_X=\Ocal_X(\kappa -d)$.  By Theorem 6.2 of
[W1], $\kappa >d$ implies the Gaussian of $C$ is surjective, so the
canonical embedding of $C$ has only conical extensions.
  \end{numpar}

  \begin{lemma} \label{lem3.8} $X$ computes the Gaussian of $C$,
except in the following cases:
   \begin{enumerate}
 \item \label{3.8.1} \qquad $n=2, \quad d_1\leq 6$
 \item \label{3.8.2} \qquad $n=3, \quad d_1\leq d_2 \leq 4$
 \item \label{3.8.3} \qquad $n=4, \quad (d_1,d_2,d_3) = (2,2,2),\;
(2,2,3),\; (2,3,3)$
 \item \label{3.8.4} \qquad $n=5, \quad \mbox{\rmfamily all }d_i =2.$
    \end{enumerate}
  \end{lemma}

  \begin{proof} $\Ocal_X(C)=\Ocal_X(d)$, and $K_X+C=\Ocal_X(\kappa )$
is very ample.  By (\ref{prop3.3}), it suffices to check that except
in the cases above, the Gaussian of $\Ocal_X(\kappa )$ is surjective
and
  \[
H^1(\Omega^1_X(2\kappa )\otimes\Ocal_C) = H^1(\Omega^1_X(2\kappa
-d))=0.
  \]
As in [W1], 6.6, the first condition follows from the surjectivity of
non-trivial Gaussians on $\Pbb^n$ plus surjectivity of the
composition
  \[
\Gamma (\Pbb^n, \Omega^1_{\Pbb}(2\kappa )) \to \Gamma
(X,\Omega^1_{\Pbb}(2\kappa )\otimes\Ocal_X) \to \Gamma
(X,\Omega^1_X(2\kappa )).
  \]
  \end{proof}

  \begin{remark} If $n\geq 2$, in all above examples $X$ is a
  canonically trivial surface (a $K$-3, if smooth), but it is not unique;
  thus $C$ could not have Gaussian of
$\corank \; 1=h^0(-K_X)$.
  \end{remark}

  \begin{theorem} \label{th3.10} Let $C\subset \Pbb^n$ be a smooth
complete intersection curve.  Aside from \ref{3.8.1}--\ref{3.8.4},
the only $C$ which could lie on a $K$-3 surface are the ``obvious
ones,'' with multi-indices
    \begin{enumerate}
 \item \label{3.10.1} \qquad $n=3, \; (4,d_2), d_2 > 4$
 \item \label{3.10.2} \qquad $n=4, \; (2,3,d_3), d_3 > 3$
 \item \label{3.10.3} \qquad $n=5, \; (2,2,2,d_4), d_4>2.$
    \end{enumerate}
In these cases, the surface $X$ is a $K$-3 (if smooth) and is a
universal extension of the canonical curve.
  \end{theorem}

  \begin{proof}  \ref{3.10.1}--\ref{3.10.3} are exactly the cases
with $\kappa =d>d_{n-2}$; the other $\kappa =d$ examples are covered
by (\ref{lem3.8}).  So Theorem \ref{th3.5} applied to $X$ gives the
result, once we show the special curves of \ref{3.1.2} are included
in (\ref{lem3.8}).  For instance, it suffices to check that the
Gaussians $\Phi (K,K^{\otimes i})$ on $C$ are surjective, all $i\geq
2$.  Restricting from $\Pbb^n$ as in (\ref{lem3.8}), one needs that
the composition
  \[
\Gamma (\Pbb^n,\Omega^1_{\Pbb^n}((i+1)\kappa )) \to \Gamma
(C,\Omega^1_{\Pbb^n}((i+1)\kappa )\otimes\Ocal_C) \to \Gamma
(C,\Omega^1_C((i+1)\kappa ))
  \]
is surjective for $i\geq 2$; a standard computation gives this
once $2\kappa >d$.  But the cases with $2\kappa\leq d$ are
included in (\ref{lem3.8}).
  \end{proof}

  \begin{remarks}
     \item Taking double covers, smooth plane curves of degree 4, 5,
or 6 lie on a $K$-3 (possibly with rational double points) --- cf.
the proof of (\ref{cor4.5}) below.
     \item Green-Lazarsfeld [GL] already had shown that a plane curve of degree
$d\geq 7$ does not lie on a $K$-3; they argue that the $g^2_d$ would
be induced by a divisor on the $K$-3 , from which a contradiction
follows.
  \end{remarks}

  \begin{proposition} \label{prop3.12} Let $C$ be a smooth curve on
the rational ruled surface $\Fbb_n$, linearly equivalent to $pB+qF$
(where $B$ is the section with $B\cdot B=-n$, $F$ is a fibre).
Assume that $p\geq 5$, and
    \begin{enumerate}
 \item[(a)] \quad if $n=0$, then $q\geq 5$
 \item[(b)] \quad if $n=1$, then $q\geq p+4$
 \item[(c)] \quad if $n=2$, then $q\geq 2p+2$.
     \end{enumerate}
Then $\Fbb_n$ computes the Gaussian of $C$, with $h^0(\Fbb_n,\Ocal
(-K)) =n+6$ if $n\geq 3$ (and $=9$ for $n\leq 3$).  Every extension
of the canonical embedding of $C$ arises from blowing up $\Fbb_n$ and
blowing down an anti-canonical divisor.  In particular, $C$ does not
sit on a $K$-3 surface.
  \end{proposition}

  \begin{proof}  The first assertion is Theorem 5.5 of [DM].  Since
$K=-2B-(n+2)F$, and $rB+sF$ is very ample iff $r\geq 1$ and $s>rn$,
we conclude that $K+C$ is very ample.  Theorem \ref{th3.5} applies
once we note $C$ can be neither trigonal, nor on a del Pezzo, nor
bielliptic,
as seen by comparing the genera and corank of the Gaussian of such
curves from [Br], [CM] with the formulas above when $C\subset\Fbb_n$.
  \end{proof}

  \begin{remarks}
 \item Blowing up $\Pbb^2$ at a point of a smooth curve $C$, the
proper transform on $\Fbb_1$ is very ample but with Gaussian of
corank 10, not 9---thus the restrictions above for small $n$.
  \item When $n\geq 3$, $-K$ has $B$ as base curve.  One sees that a
generic extension of $C\subset \Pbb^{g-1}$ has a non-smoothable cusp
singularity with two exceptional curves in the minimal resolution.
  \end{remarks}

  \begin{numpar}  Now suppose $C$ is bielliptic, of genus $\geq 6$.
Then $C\subset\Pbb^{g-1}$ is the complete intersection of a quadric
$Q$ with $X$, a projective cone over an elliptic curve $D$ of degree
$g-1$ in $\Pbb^{g-2}$ [CM]. It is known that
   \begin{alignat*}{2}
&\dim H^0(C,N_{C/\Pbb}(-1))/\Gamma (\Theta_{\Pbb}(-1)) &&= 2g-2 \\
&\dim H^0(C,N_{C/\Pbb}(-i)) &&= 1 \qquad\qquad i=2 \\
& &&=0 \qquad\qquad i\geq 3.
    \end{alignat*}
$D\subset \Pbb^{g-2}$ satisfies $(N_2)$; it is defined by quadratic
equations $f_1,\cdots ,f_N$, in variables $x_2,\cdots ,x_g$.  We may
assume coordinates chosen in $\Pbb^{g-1}$ so that $C$ is defined by
these equations plus another of the form $h=x^2_1 - A(x_2,\cdots
,x_g)$, where $A$ is also quadratic.
  \end{numpar}

  \begin{theorem} \label{th3.15} (cf. [Re]).  Suppose $C$ is
bielliptic of genus $g\geq 11$.  Then there is a unique nontrivial
extension of the canonical embedding of $C$, a birationally ruled
surface with two simple elliptic singularities of degree $(g-1)$.  In
particular, $C$ does not sit on a $K$-3 surface.
  \end{theorem}

  \begin{proof}  $C$ is defined by $f_1,\cdots ,f_N$, and $h$, and
the relations are generated by the linear relations among the $f_i$,
plus the trivial relations $f_i\cdot h-h\cdot f_i =0$ (since $h$ is
not a 0-divisor mod the $f_i$).  Let $W\subset\Pbb^g$ be an extension
of $C\subset\Pbb^{g-1}$, defined by $t=0$.  Then liftings of the
equations $f_i$ induce an extension of $X$ as well.  But it is
well-known ([P]) that the cone over an elliptic curve of degree
$\geq 10$ has no non-trivial extensions; thus, when $g-1 \geq 10$,
after changing coordinates one may assume the equations defining $W$
are given by $f_1,\cdots f_N, h+t\alpha +t^2\beta$.  Given the form
of $h$, replacing $x_1$ by $x_1-t\alpha /2$ changes only the last
equation, to $h+t^2\beta$, where $\beta\in\Cbb$.  If $\beta =0$, one
has a cone; otherwise, replacing $t$ by $t/\sqrt{\beta}$, one has the
unique non-trivial extension of $C$ (with $\beta =1$).  This is
described geometrically via Corollary \ref{cor3.7}.  Noting $K_X
\cong \Ocal_X(1)$, let $E=C$ and $F\cong D$ a smooth hyperplane
section of $X\subset\Pbb^{g-1}$.  Blowing up the $2(g-1)$ points of
$C\cap D$, then blowing down the proper transform of $D$, gives an
extension of $C$ with a second elliptic singularity, also of degree
$g-1$.  Note that the $KS$ map is 0 for this extension.
  \end{proof}

       \section{Complete intersections which lie on Calabi-Yaus}

  \begin{numpar} Let $X\subset\Pbb^n$ be a smooth complete
intersection subvariety of dimension $r\geq 2$, of multidegree $2\leq
d_1\leq \cdots \leq d_{n-r}$.  Assume $\kappa =\sum d_i-(n+1) >0$, so
that $X$ has very ample canonical bundle $\Ocal_X(\kappa )$.  Let
$Y_{\alpha}\supset X$ be a general complete intersection defined by
the same equations as $X$ except for one of degree $d_{\alpha}$.
Then $K_{Y_{\alpha}} \cong \Ocal_{Y_{\alpha}} (\kappa -d_{\alpha})$.
Thus, if $\kappa =d_{\alpha}$ for some $\alpha$, then $X$ is a
hyperplane section of a Calabi-Yau, assuming $Y_{\alpha}$ can be
chosen to have
rational singularities.  More generally,
  \end{numpar}

  \begin{proposition} \label{prop4.2} Let $X$ be a complete
intersection so that $\kappa$ divides $d_{\alpha}$, for some
$\alpha$.  Suppose $Y=Y_{\alpha}$ can be chosen to have at most
rational singularities. Then there is a cyclic branched cover of $Y$
which is a Calabi-Yau variety and on which $X$ is a hyperplane
section.
  \end{proposition}

  \begin{proof} If $m\kappa =d_{\alpha}$, then
$\Ocal_Y(X)=\Ocal_Y(m\kappa )$; form the $m$-fold cyclic branched
cover of $Y$ along $X$.  It is easy to check $Y$ is a Calabi-Yau
variety, and the reduced ramification divisor $X'\cong X$ is very
ample.
  \end{proof}

  \begin{numpar} If $\kappa$ is large compared to the $d_{\alpha}$,
then $X$ cannot sit on a Calabi-Yau, as is proved for curves in
(\ref{th3.10}).  Specifically,
  \end{numpar}

  \begin{theorem} \label{th4.4}  Let $X\subset\Pbb^n$ be a smooth
complete intersection, with multidegree $d_1\leq \cdots \leq
d_{n-r-1} \leq d_{n-r}=d$, and with $\kappa =\sum d_i-(n+1) > 0$.
    \begin{enumerate}
 \item[(a)] If $\kappa >d$, then $(X,K_X)$ has no extensions.
 \item[(b)] If $d>\kappa >\frac{d}{2}$ and $\kappa >d_{n-r-1}$, then
$(X,K_X)$ has non-trivial extensions, but all have non-rational
singularities.
     \end{enumerate}
  \end{theorem}

  \begin{corollary} \label{cor4.5}  A smooth hypersurface $X\subset
\Pbb^n$ of degree $d$ is a hyperplane section of a Calabi-Yau variety
iff $n+1 \leq  d\leq 2n+2$.
  \end{corollary}

  \begin{corollary} \label{cor4.6} Let $X\subset\Pbb^n$ be a smooth
complete intersection of type $(d' ,d)$, with $d'\leq d$.  Then
    \begin{enumerate}
 \item[(a)] if $d\leq n+1$, $X$ sits on a smooth Calabi-Yau
hypersurface.
 \item[(b)] if $d>n+2$ and $\frac{d}{2} +d' >n+1$, then $X$ cannot
sit on a Calabi-Yau variety unless $d' =n+1$.
 \item[(c)] if $\frac{d}{2} +d' =n+1$, then the generic complete
intersection variety of type $(d',d)$ sits on a Calabi-Yau.
    \end{enumerate}
  \end{corollary}

  \begin{proof}[Proof of 5.4] We assume $2\kappa >d$.  Standard
  cohomological arguments yield,
  \begin{align*}
H^1(\Theta_X(-i\kappa )) &=0, \vspace{2.3in} i\geq 2,  \\
H^1(\Theta_X(-\kappa )) &\cong \oplus \Gamma (\Ocal_X(d_i-\kappa )).
  \end{align*}
It follows from the normal bundle sequence that the canonical
embedding of $X$ satisfies $\Gamma (N(-2))=0$.  (a) follows easily
from (1.9) and (2.3.2). Further, $\kappa>d_{n-r-1}$ implies
  \[
H^1(\Theta_X(-\kappa )) \cong \Gamma (\Ocal_X(d-\kappa )).
  \]
With $Y=Y_{n-r}$ as above, this last term equals $h^0(Y,-K_Y)$, and
one easily finds that $H^0(\Theta_Y(-(K_Y+X))\otimes\Ocal_X)=0$.  Now
apply Corollary \ref{cor3.12}.
  \end{proof}

  \begin{proof}[Proof of 5.5] If $d=2n+2$, Proposition 4.2 applies;
so it remains to show that $X$ is on a
Calabi-Yau variety when $n+1 < d < 2n+2$.  Let $X'$ be a smooth
hypersurface of degree $d' =2n+2-d$ intersecting $X$ transversally.
Let $Z\to \Pbb^n$ be the blow-up along $X\cap X'$, with $\bar{X}$ and
$\bar{X}'$ the proper transforms; take the double cover $Y\to Z$
branched along $\bar{X}\cup \bar{X}'$.  Then $Y$ is a Calabi-Yau
variety, with a divisor $\tilde{X}$ isomorphic to $X$, for which the
normal bundle is $K_X$.  In particular, $\Ocal_Y(\tilde{X})$ is very
ample in a neighborhood of $\tilde{X}$, and defines a morphism $Y\to
Y^*$ which collapses $\tilde{X}'$ to a rational singular point (since
$d' < n+1$).  We can view $X$ as a hyperplane section of the singular
Calabi-Yau $Y^*$ (cf. Lemma 3.4).
  \end{proof}

  \begin{proof}[Proof of 5.6] If $d\leq n+1$, by Bertini's Theorem
the general hypersurface of degree $n+1$ containing $X$ is smooth,
and is a Calabi-Yau.  (b) follows easily from (\ref{th4.4}), since
$\kappa =d' +d-(n+1)$.  The condition in (c) may be written $2\kappa
=d$, so (\ref{prop4.2}) applies.
  \end{proof}

  \end{document}